\documentclass[12pt,leqno]{amsart}  
\usepackage{pgf,tikz}
\usepackage{amsmath,amstext,amsthm,amssymb,amsxtra}
\usepackage{txfonts} 
\usepackage[pagebackref,hypertexnames=false]{hyperref} 
\usepackage[backrefs]{amsrefs}

\allowdisplaybreaks
\swapnumbers

\allowdisplaybreaks

\theoremstyle{plain} 
 
\newtheorem{proposition}[equation]{Proposition} 
\newtheorem{theorem}[equation]{Theorem} 
 
\newtheorem{nehari}[equation]{Nehari's Theorem}

\theoremstyle{definition}

\theoremstyle{remark}

\numberwithin{equation}{section}

\def\norm#1.#2.{\lVert#1\rVert_{#2}}
\def\Norm#1.#2.{\bigl\lVert#1\bigr\rVert_{#2}}
\def\NOrm#1.#2.{\Bigl\lVert#1\Bigr\rVert_{#2}}
\def\NORm#1.#2.{\biggl\lVert#1\biggr\rVert_{#2}}
\def\NORM#1.#2.{\Biggl\lVert#1\Biggr\rVert_{#2}}

\def\ip#1,#2,{\langle #1,#2\rangle}
\def\Ip#1,#2,{\bigl\langle#1,#2\bigr\rangle}
\def\IP#1,#2,{\Bigl\langle#1,#2\Bigr\rangle}

\def\mid{\,:\,}

\def\ABs#1{\biggl\lvert#1\biggr\rvert}

\def\XXint#1#2#3{{\setbox0=\hbox{$#1{#2#3}{\int}$}
     \vcenter{\hbox{$#2#3$}}\kern-.5\wd0}}

\def\eqdef{\stackrel{\mathrm{def}}{{}={}}}

\title {Haar Shifts, Commutators,  and Hankel Operators}
\author{Michael Lacey}   

\address{School of Mathematics \\ Georgia Institute of Technology \\ Atlanta GA 30332 }

\email{lacey@math.gatech.edu}

\begin{document}

\begin{abstract}
{Hankel operators lie at the junction of analytic and real-variables.  We will explore this junction, 
from the point of view of Haar shifts and commutators. }
\end{abstract} 

\maketitle

\section{Haar Functions} 

We consider operators which satisfy  invariance properties with respect to two 
well-known groups. The first group 
we take to the \emph{translation} operators 
\begin{equation}\label{e.trans}
\operatorname {Tr}_y f (x) := f (x-y)\,, \qquad y\in \mathbb R \,. 
\end{equation}
Note that formally, the adjoint operator is $ (\operatorname {Tr}_y) ^{\ast} = 
\operatorname {Tr} _{-y}$.   The collection of operators $ \{\operatorname {Tr}_y \mid y\in \mathbb R \}$ 
is a representation of the additive group $ (\mathbb R , +)$.

It is an important, and very general principle that a linear operator $ \operatorname L$ acting 
on some vector space of functions, which is assumed to commute with all translation operators, is 
in fact given as convolution, in general with respect to a measure or distribution, thus, 
\begin{equation*}
\operatorname L f (x) = \int f (x-y) \; \mu (dy)\,. 
\end{equation*}
For instance, with the identity operator, $ \mu $ would be the Dirac pointmass at the origin.  

The second group is the set of \emph{dilations on $ L ^{p}$}, given by 
\begin{equation}\label{e.dil}
\operatorname {Dil} ^{ (p)} _{\lambda } f (x) :=  \lambda ^{-1/p} f (x/\lambda )\,, 
\qquad 
0< \lambda , p<\infty \,. 
\end{equation}
Here, we make the definition so that $ \norm f. p. = \norm \operatorname {Dil} _ {\lambda } ^{(p)} f .p. $. 
  The \emph{scale} of the dilation $ \operatorname {Dil} _ {\lambda
} ^{(p)} $  is said to be $ \lambda $, and these operators are a representation of the
 multiplicative group $ (\mathbb R _+ , \ast  )$. 
The Haar measure of of this group is $ dy/y$.

Underlying this subject are the delicate interplay between local averages and differences.  Some of this 
interplay can be encoded into the combinatorics of \emph{grids}, especially the \emph{dyadic grid,} defined 
to be  $\mathcal D:=\{2^{k}(j,j+1)\mid j,k\in\mathbb Z\}$.   

The Haar functions are a remarkable class of functions indexed by the dyadic grid $ \mathcal D$.  Set 
\begin{equation*}
h (x)= -\mathbf 1_{(-1/2,0)}+\mathbf 1_{(0,1/2)}\,, 
\end{equation*}
a mean zero function supported on the interval $ (-1/2,1/2)$, taking two values, with $ L ^{2}$ norm equal to one. 
Define the \emph{Haar function } (associated to interval $ I$) to be 
\begin{align}\label{e.haardef}
h _I &:= \operatorname {Dil} _{I} ^{2} h_I 
\\ \label{e.dilIdef}
\operatorname {Dil}_I ^{(p)} &:= \operatorname {Tr}_{c (I)}  \operatorname {Dil} _{\lvert  I\rvert } ^{ (p)} 
\,, \qquad c (I)= \textup{center of $ I$.}
\end{align}
Here, we introduce the notion for the \emph{Dilation associated with interval $ I$.}

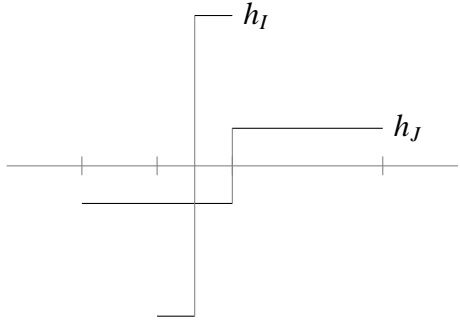
\begin{figure}
 \begin{tikzpicture}
 \draw[gray] (-1,0) -- (5,0); 
 \draw[gray] (0,-.125) -- (0,.125);
 \draw[gray] (4,-.125) -- (4,.125);
 \draw[gray] (1,-.125) -- (1,.125);
 \draw[gray] (2,-.125) -- (2,.125);
 \draw  (0,-.5) -- (2,-.5);
 \draw[gray] (2,-.5) -- (2,.5);
 \draw  (2,.5) -- (4,.5) node[right] {$ h_J$}; 
 \draw  (1,-2) -- (1.5,-2); 
 \draw[gray] (1.5,-2) -- (1.5,2); 
 \draw  (1.5,2) -- (2,2) node[right] {$ h_I$};
 \end{tikzpicture}
\caption{Two Haar functions.}
\label{f.haar}
\end{figure}

The Haar functions have profound properties, due to their connection to both analytical and probabilistic 
properties.   An elemental property is that they form a basis for $ L ^2 (\mathbb R )$. 

\begin{theorem}\label{t.haarBasis} The set of functions $ \{\mathbf 1_{[0,1]}\}\cup 
\{ h_I \mid I \in \mathcal D, \ I\subset [0,1]\}$ form an orthonormal basis for $ L ^{2} ([0,1])$.  
The set of functions $ \{h_I \mid I\in \mathcal D\}$ form an orthonormal basis for $ L ^{2} (\mathbb R )$. 

\end{theorem}

\section{Paraproducts} 

Products, and certain kind of renormalized products are common objects.  Let us explain 
the renormalized products in a very simple situation.  We begin with the definition of 
a \emph{paraproduct}, as a bilinear operator.  Define 
\begin{align}\label{e.h01}
h _{I} ^{0}= h _I\,, \qquad h_I ^{1}= \lvert  h_I ^{0}\rvert = \operatorname {Dil} _{I} ^{2} \mathbf 1_{[-1/2,1/2]} \,.  
\end{align}
The superscript $ {}^{0}$ indicates a mean-zero function, while the superscript $ {} ^{1}$ indicates a 
non-zero integral.   Now define 
\begin{equation}\label{e.paraproductDefinition}
\operatorname P ^{\epsilon _1, \epsilon _2, \epsilon _3} (f_1,f_2) 
\coloneqq 
\sum _{I\in \mathcal D} \frac {\ip f_1, h _{I} ^{\epsilon _1}, } {\sqrt {\lvert  I\rvert }}
\ip f_2, h_2 ^{\epsilon _2} , h _{I} ^{\epsilon _3} \,, \qquad \epsilon _j\in \{0,1\}.  
\end{equation}
For the most part, we consider cases where there is one choice of $ \epsilon _j$ which is equal to one, 
but in considering fractional integrals, one considers examples where all $ \epsilon _j$ are  equal to one. 
The triple $ (\epsilon _1, \epsilon _2, \epsilon _3) $ is the \emph{signature} of the Paraproduct. 

We have chosen this definition for specificity, but at the same time, it must be stressed that there is no 
canonical definition, and the presentation of a paraproduct can differ in a number of ways.  Whatever the presentation, 
its single most important attribute is its signature.  Indeed, in Proposition~\ref{p.[]=}, we will see that 
a paraproduct arises from a computation that, while not of the form above, is clearly an operator of 
signature $ (0,0,0)$.  All the important prior work on commutators,  see \cites{MR1349230
,MR631751
,MR511821
,MR518170
,MR54:843}
can be interpreted in this notation.   (The Lectures of M.~Christ \cite{MR1104656} are recommended as 
a guide to this literature.)    
For instance, in the notation of Coifman and Meyer \cites{MR518170,MR511821}, a $ P_t$ denotes a $ {}^1$, 
while a $ Q_t$ denotes a $ {}^0$.

Why the name paraproduct? This is probably best explained by the identity 
\begin{equation}\label{e.product=3para}
f_1 \cdot f_2 = 
\operatorname P ^{1,0,0} (f_1,f_2)
+
\operatorname P ^{0,0,1} (f_1,f_2)+\operatorname P ^{0,1,0} (f_1,f_2) \,. 
\end{equation}
Thus, a product of two functions is a sum of three paraproducts.  The three individual paraproducts 
in many respects behave like products, for instance we will see that there is a H\"older Inequality.  
And, very importantly, in certain instances they are \emph{better} than a product.  

To verify \eqref{e.product=3para}, let us first make the self-evident observation that 
\begin{equation}\label{e.avgInTermsOfHaar}
\frac 1 {\lvert  J\rvert }\int _J g (y)\; dy 
= \frac {\ip g, h ^{1}_I, } {\sqrt {\lvert  I\rvert }} 
= \sum _{J \;:\; J\supsetneq I} \ip g, h_J,  h_J (I)\,, 
\end{equation}
where $ h_J (I)$ is the (unique) value $ h_J$ takes on $ I$.  In \eqref{e.product=3para}, 
expand both $ f_1$ and $ f_2$ in the Haar basis,
\begin{align*}
 f_1 \cdot f_2  & = 
\Biggl\{ \sum _{I\in \mathcal D} \ip f_1, h_I, h_I \Biggr\}
\cdot 
\Biggl\{ \sum _{J\in \mathcal D} \ip f_2, h_J, h_J \Biggr\}\,. 
\end{align*}
Split the resulting product into three sums, (1) $I=J $, (2) $ I\subsetneq J$ (3) $ J\subsetneq I$. 
In the first case, 
\begin{equation*}
\sum _{I,J \;:\; I=J}  \ip f_1, h_I,\ip f_2, h_J, (h_I) ^2= 
\operatorname P ^{0,0,1} (f_1,f_2) \,. 
\end{equation*}
In the second case, use \eqref{e.avgInTermsOfHaar}. 
\begin{align*}
\sum _{I,J \;:\; I\subsetneq J}  \ip f_1, h_I,\ip f_2, h_J, h_I  \cdot  \tfrac 1 {\lvert  I\rvert }  
\int _I 
h_J (y) \; dy  
&
=\sum _{I} \ip f_1, h_I,   \frac { \ip f_2, h _I ^{1}, } {\sqrt {\lvert  I\rvert }} h_I 
\\
{= \operatorname P ^{0,1,0} (f_1,f_2)} \,. 
\end{align*}
And the third case is as in the second case, with the role of $ f_1$ and $ f_2$ switched. 

A rudimentary property is that Paraproducts should respect H\"older's inequality, a matter that we turn to next. 
This Theorem is due to Coifman and Meyer \cites{MR511821,MR518170}.  Also see \cites{camil1,camil2,math.CA/0502334}.

\begin{theorem}\label{t.paraproductHolder}  Suppose at most one of $ \epsilon _1, \epsilon _2, \epsilon _3$ 
are equal to one. We have the inequalities 
\begin{equation}\label{e.paraproductHolder} 
\norm \operatorname P ^{\epsilon _1, \epsilon _2, \epsilon _3} (f_1, f_2) . q. 
\lesssim 
\norm f_1. p_1. \norm f_2. p_2 .  \,, 
\qquad 
1<p_1,p_2<\infty \,,\, 1/q=1/p_1+1/p_2 \,. 
\end{equation}

\end{theorem}
\section{Paraproducts and Carleson Embedding} 

We have indicated that  Paraproducts are better than products in one way.  These fundamental  inequalities
are the subject of this section.  Let us define the notion of \emph{(dyadic) Bounded Mean Oscillation}, $ \textup{BMO}$
for short, by 
\begin{equation}\label{e.BMOdef} 
\norm f. \textup{BMO}. = \sup _{J\in \mathcal D} \Biggl[ \lvert  J\rvert ^{-1} \sum _{I\subset J} 
\ip f, h_R, ^2 \Biggr] ^{1/2} \,. 
\end{equation}

\begin{theorem}\label{t.ParaproductBMO} Suppose that at exactly one of $ \epsilon _2$ and $ \epsilon _3$ 
are equal to $ 1$.  
\begin{equation}\label{e.ParaproductBMO} 
\Norm \operatorname P ^{0,\epsilon _2,\epsilon _3} (f_1, \cdot ) . p\to p.
\simeq \norm f_1. \textup{BMO}.  \,, \qquad 1<p<\infty \,. 
\end{equation}
Indeed, we have 
\begin{equation}
\Norm \operatorname P ^{0,1,0} (f_1, \cdot ) . p\to p.
\simeq \sup _J \Norm \operatorname P ^{0,1,0} (f_1, \lvert  J\rvert ^{-1/p} \mathbf 1_{J} ) . p. 
\simeq \norm f_1. \textup{BMO}. \,. 
\end{equation}
\end{theorem}

Here, we are treating the paraproduct as a linear operator on $ f_2$, and showing that the 
operator norm is characterized by $ \norm f_1. \textup{BMO}. $. 
Obviously, $ \norm f. \textup{BMO}. \le 2 \norm f. \infty .$, and again this a crucial point, 
there are unbounded functions with bounded mean oscillation, with the canonical example being $ \ln x$. 
Thus, these paraproducts are, in a specific sense, better than pointwise products of functions.

\begin{proof}
The case $ p=2$ is essential,  and the only case considered in these notes.  This particular case is 
frequently referred to as \emph{Carleson Embedding}, a term that arises from the original 
application of the principal in the Corona Theorem. 

Let us discuss the case of $ \operatorname P ^{0,1,0}$ in detail.  
Note that the dual of the operator 
\begin{equation*}
f_2 \longrightarrow \operatorname P ^{0,1,0} (f_1, f_2)\,, 
\end{equation*}
that is we keep $ f_1$ fixed, is the operator $ \operatorname P ^{0,0,1} (f_1, \cdot )$, so it is enough 
to consider $ \operatorname P ^{0,1,0} $ in the $ L ^{2}$ case. 

One direction of the inequalities is as follows.   
\begin{align*}
\norm \operatorname P ^{0,\epsilon _2,\epsilon _3} (f_1, \cdot ) . 2\to 2.
& \ge \sup _{J} \norm \operatorname P ^{0,\epsilon _2,\epsilon _3} (f_1, h ^{1}_J ) . p.
\\
& \ge \norm f_1. \textup{BMO}. 
\end{align*}
as is easy to see from inspection.  Thus, the $ \textup{BMO}$ lower bound on the operator norm arises 
solely from testing against normalized indicator sets. 

For the reverse inequality, 
we compare to the Maximal Function. Fix $ f_1, f_2$, and let 
\begin{gather*}
{\mathcal D_k = 
\{ I\in \mathcal D \mid   \frac { \lvert\ip f_2, h_I,\rvert } { \sqrt{\lvert  I\rvert }} \simeq 2 ^{k}\} }
\end{gather*}
Let $ \mathcal D_k ^{{\ast}} $ be the maximal intervals in $ \mathcal D_k$. The $ L ^{2}$-bound for the Maximal 
Function gives us 
\begin{gather} \label{e.tmf}
\sum _{k} 2 ^{2k } \sum _{ I ^{{\ast}} \in \mathcal D_k ^{{\ast}} } \lvert  I ^{\ast} \rvert \lesssim  
\norm \operatorname M f_2. 2. ^2 \lesssim \norm f.2. ^2 \,. 
\end{gather}
 Then, for $ I ^{{\ast}} \in \mathcal D _k ^{{\ast}} $
we have 
\begin{align*}
\NOrm \sum _{\substack{I\in \mathcal D_k\\ I\subset I ^{{\ast}}  }} 
\ip f _1, h_I,   2 ^{k} h_I .2. ^2 
&=  
 2 ^{2k}\sum _{\substack{I\in \mathcal D_k\\ I\subset I ^{{\ast}}  }   } \ip f _1, h_I, ^2  
\\
& \le 2 ^{2k} \norm f_1. \textup{BMO}. ^2  \lvert  I ^{{\ast}} \rvert 
\end{align*}
And so we are done by \eqref{e.tmf}. 

\end{proof}

\section{Hilbert Transform} 

 It is a useful Theorem, one that we shall return to later, that the set of operators $ \operatorname L$ that 
 are bounded from $ L ^{2} (\mathbb R )$ to itself, and commute with both translations and dilations have a 
 special form.  They are linear combinations of the Identity operator, and the \emph{Hilbert } transform.  
 The latter operator, fundamental to this study, is given by 
 \begin{equation} \label{e.hilbertDef}
\operatorname H f (x) := \textup{p.v.} \int f (x-y) \; \frac {dy} y \,. 
\end{equation}
Here, we take the integral in the \emph{principal value} sense, as the kernel $ 1/y$ is not integrable.  
Taking advantage of the fact that the kernel is odd, one can see that the limit below 
\begin{equation} \label{e.pvDef}
\lim _{\epsilon \to 0} \int _{\epsilon < \lvert  y\rvert < 1/ \epsilon  } f (x-y) \; \frac {dy}y 
\end{equation}
exists for all $ x$, provided $ f$ is a Schwartz function, say.   Thus, $ \operatorname H$ has 
an unambiguous definition on a dense class of functions, in all $ L ^{p}$.    
We shall take \eqref{e.pvDef} as our general definition of principal value.  The Hilbert transform is the 
canonical example of a \emph{singular integral}, that is one that has to be defined in some principal value sense.

Observe that $ \operatorname H$, being convolution commutes with all translations.  That is also commutes with 
all dilation operators follows from the observation that $ 1/y$ is a multiple of the multiplicative Haar measure.  
It can also be recovered in a remarkably transparent way from a simple to define 
operator based upon the Haar functions.   Let us define 
\begin{align}\label{e.gIdef}
g &= - \mathbf 1_{(-1/4,-1/4)}+ \mathbf 1_{(-1/4,1/4)} - \mathbf 1_{(1/4,1/2)}
\\& = 
2 ^{-1/2} \{ h _{(-1/2,0)} + h _{(0,1/2)} \} 
\\
\mathfrak H f &= \sum _{I\in \mathcal D} \ip f, h_I, g_I \,, 
\end{align}
where as before, $ g _I  = \operatorname {Dil} _{I} ^{(2)}g$.   It is clear that $ \mathfrak H$ is a bounded 
operator on $ L ^2 $.  What is surprising is that that it can be used to recover the Hilbert transform 
exactly.   The succinct motivation for this definition is that $ \operatorname H (\sin )= \cos$, so that 
if $ h_I$ is a local sine, then $ g_I$ is a local cosine.

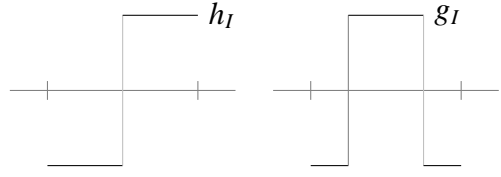
\begin{figure}
\begin{center}
\begin{tikzpicture}
\draw[gray] (-1.5,0) -- (1.5,0); 
\draw[gray] (-1,-.125) -- (-1,.125);
\draw[gray] (1,-.125) -- (1,.125);
\draw (-1,-1) -- (0,-1) ;
\draw[lightgray] (0,-1) -- (0,1);
\draw (0,1) -- (1,1)  node[right] {$ h_I$} ;
\draw[gray] (2,0) -- (5,0); 
\draw[gray] (2.5,-.125) -- (2.5,.125);
\draw[gray] (4.5,-.125) -- (4.5,.125);
\draw (2.5,-1) -- (3,-1) ;
\draw[gray] (3,-1) -- (3,1);
\draw (3,1) -- (4,1) node[right] {$ g_I$} ; 
\draw[lightgray] (4,1) -- (4,-1);
\draw (4,-1) -- (4.5,-1); 
\end{tikzpicture}
\end{center}
\caption{A Haar function $ h_I$ and its dual $ g_I$.}
\label{f.h&g}
\end{figure}

\begin{theorem}[S.~Petermichl \cite{MR1756958}]
\label{t.avergeHilbert}  There is a non-zero constant $ c$ so that 
\begin{equation} \label{e.hilberShift}
H = c \lim _{Y\to \infty } \int _0 ^{Y} \int _1 ^{2} 
\operatorname {Tr} _{y} \operatorname {Dil} _{\lambda } ^{(2)} 
\mathfrak H 
 \operatorname {Dil} _{1/\lambda } ^{(2)}  \operatorname {Tr} _{-y} 
\;  \frac {d \lambda } {\lambda }\, \frac {dy} Y \,. 
\end{equation}
\end{theorem}

As a Corollary, we have the estimate $ \norm H .2. \lesssim 1$, as $  \mathfrak H$ is clearly bounded on $ L ^{2}$. 

The operator $ \mathfrak h$ is referred to as a \emph{Haar shift} or as a \emph{dyadic shift} (\cite{MR1964822}).
Certain canonical singular integrals, like the Hilbert, Riesz and Beurling transform admit remarkably
simple Haar shift variants, which fact can be used to prove a range of deep results.  See for instance \cites{MR2367098
,MR2354322
,MR1894362
,MR2164413}.  For applications of this notion to more general singular integrals, see \cite{onBeyondCRW}*{Section 4}.

\begin{proof}
Consider the limit on the right in (\ref{e.hilberShift}). 
This is seen to exist for each $ x\in \mathbb R $ for Schwartz functions $ f$.
While this is elementary, it might be useful for us to define the 
auxiliary operators 
\begin{equation*}
\operatorname T_j f \coloneqq \sum _{\substack{I\in \mathcal D\\ \lvert  I\rvert \le 2 ^{j} }}
\ip f, h_j, g_j \,. 
\end{equation*}
The individual terms of this series are rapidly convergent.  As $ \lvert  I\rvert  $ 
becomes small, one uses the smoothness of the function $ f$.  As $ \lvert  I\rvert $ 
becomes large, one uses the fact that $ f$ is integrable, and decays rapidly.  
Call the limit $ \widetilde {\operatorname H}f$.  

Let us also note that the operator $ \operatorname T_j  $ is invariant under 
translations by an integer multiple of $ 2 ^{j}$. Thus, the auxiliary operator 
\begin{equation*}
2 ^{-j}\int _{0} ^{2 ^{j}} \operatorname {Tr} _{-t} f \operatorname {Tr} _{t} \; dt 
\end{equation*}
will be translation invariant. 
  Thus $  \widetilde {\operatorname H}$ is convolution with respect to 
a linear functional on Schwartz functions, namely a distribution.  

Concerning dilations, $ \operatorname T$ is invariant under dilations 
by a power of $ 2$.  
Now, dilations form a group under multiplication on $ \mathbb R _+$, and this 
group has Haar measure $ d\delta/ \delta  $ so that the operator below 
will commute with all dilations. 
\begin{equation*}
\int _{0} ^{1} 
\operatorname {Dil} _{1/\delta } ^{2}  
\operatorname T \operatorname {Dil} _{\delta } ^{2} \frac {d \delta } \delta 
\end{equation*}
Thus, $ \widetilde {\operatorname H}$ 
commutes with all dilations.  

Therefore, $ \widetilde {\operatorname H}$ must be a linear combination of 
a Dirac delta function and convolution with $ 1/y$. (The function 
$ 1/\lvert y  \rvert $ is also invariant under dilations, but the inner product 
with this function is not a linear functional on distributions.) 
Applying $ \widetilde {\operatorname H}$ to a non negative Schwartz function 
yields a function with zero mean.  Thus, $ \widetilde {\operatorname H}$ 
must be a multiple of convolution with $ 1/y$, and we only need to see that 
it is non zero multiple. 

Let us set $ G_j$ to be the operator 
\begin{equation*}
\operatorname G_j f \coloneqq   \int _{0} ^{2 ^{j}} 
\operatorname {Tran} _{t} \sum _{\substack{I\in \mathcal D\\ \lvert  I\rvert= 2 ^{j}  }}
\ip \operatorname {Tran} _{-t}f , h_I, h_I \; \frac{dt} {2 ^{j}} \,. 
\end{equation*}
This operator translates with translation and hence is convolution.  
We can write $ \operatorname G_j f= \gamma _j \ast f$.  By the dilation invariance of the 
Haar functions, we will have $ \gamma _j= \operatorname {Dil} _{2 ^{j}} ^{1} \gamma _0$.
A short calculation shows that 
\begin{equation*}
\gamma _0 (y)=\int _{0} ^{1} h_I (y+t) h_I (y)\; dt 
\end{equation*}
This function is depicted in Figure~\ref{f.hilbertshift}.  
Certainly  the operator $ \sum _{j} G_j $ is convolution with 
$ \sum _{j} \gamma _j (x)$.  This kernel is odd and is strictly positive on $ [0,\infty )$. 
This finishes our proof.

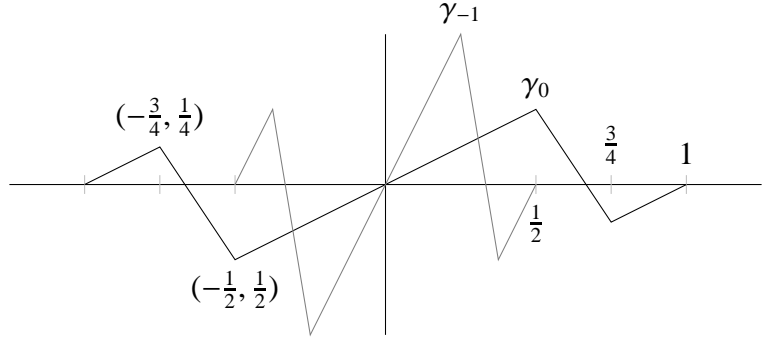
\begin{figure} 
\begin{center}
\begin{tikzpicture}
\draw (-5,0) -- (5,0) ;
\draw (0,-2) -- (0,2); 
\draw[lightgray] (-4,-.125) -- (-4,.125) ; 
\draw[lightgray] (-3,-.125) -- (-3,.125) ; 
\draw[lightgray] (-2,-.125) -- (-2,.125) ; 
\draw[lightgray] (4,-.125) -- (4,.125) node [above,black] {1} ; 
\draw[lightgray] (3,-.125) -- (3,.125) node [above,black] {$ \tfrac34$} ; 
\draw[lightgray] (2,.125) -- (2,-.125) node [below,black] {$ \tfrac12$} ; 
\draw (-4,0) -- (-3,.5) node [above] {$ (-\tfrac34,\tfrac14)$} -- (-2,-1)
node [below] {$ (-\tfrac12,\tfrac12)$}--  (2,1) node[above] {$ \gamma _0$} -- (3,-.5) -- (4,0); 
\draw[gray] (-2,0) -- (-1.5,1) -- (-1,-2) 
--  (1,2)  node[above,black] {$ \gamma _{-1}$} -- (1.5,-1) -- (2,0); 
\end{tikzpicture}
\end{center}
\caption{The graph of $ \gamma _0$ and $ \gamma _{-1}$.}
\label{f.hilbertshift}
\end{figure}

\end{proof}

\section{Commutator Bound} 

We would like to explain a classical result on commutators. 

\begin{theorem}\label{t.comm} For a function $ b$, and $ 1<p<\infty $ we have the equivalence 
\begin{equation*}
\norm [b, \operatorname H] . p\to p. \simeq 
\norm b. \textup{BMO}. \,, 
\end{equation*}
where this is the non-dyadic $ \textup{BMO}$ given by 
\begin{equation*}
\sup _{\textup{$ I$ interval}} \Biggl[ \lvert  I\rvert ^{-1} \int _{I} 
\ABs{f - \lvert  I\rvert ^{-1} \int_I f (y)\; dy } \; dx\Biggr] ^{1/2} \,. 
\end{equation*}
\end{theorem}

We refer to this as a classical result, as it can be derived from the Nehari theorem, as we will explain below.  
The lower bound on the operator norm is found by applying the commutator to normalized indicators of integrals, 
and we suppress the proof. 

Both  bounds are   very easy, if one appeals to the Nehari Theorem.  See our comments on Nehari's Theorem below. 
But, in many circumstances, different proofs admit different modifications, 
and so we present a `real-variable' proof,  deriving the upper bound from the Haar shift, and 
the Paraproduct bound in a transparent way.  

Replacing the Hilbert transform by the Haar Shift, we prove 
\begin{equation}\label{e.shift-comm}
\norm [ b, \mathfrak H] . p\to p . \lesssim \norm b. \textup{BMO}. 
\end{equation}
The last norm is dyadic-$ \textup{BMO}$, which is strictly smaller than non-dyadic $ \textup{BMO}$.
But Theorem~\ref{t.avergeHilbert} requires that we use all translates and 
dilates to recover the Hilbert transform, and so the non-dyadic $ \textup{BMO}$ norm will be invariant under these 
translations and dilations.  

The Proposition is that $  [ b, \mathfrak H]$ can be explicitly computed as a sum of Paraproducts which 
are bounded. 

\begin{proposition}\label{p.[]=}  We have 
\begin{align}  \label{e.px}
 [ b, \mathfrak H]f&= \operatorname P ^{0,1,0} (b, \mathfrak H f)- \mathfrak H \circ \operatorname P ^{0,1,0} (b, f) 
\\
& \qquad + \operatorname P ^{0,0,1} (b, \mathfrak H f)- \mathfrak H \circ \operatorname P ^{0,0,1} (b, f)  
\\ \label{e.oddONe}
& \qquad + \widetilde {\operatorname  P} ^{0,0,0} (b, f) \,. 
\end{align}  
In the last line, $ \widetilde {\operatorname  P} ^{0,0,0} (b, f)$ is defined to be 
\begin{equation*}
\widetilde {\operatorname  P} ^{0,0,0} (b, f)
= 
\sum _{I\in \mathcal D} \frac {\ip b, h _{I} ^{0},} {\sqrt I} \ip f, h_I ^{0},  ( h_{I _{\textup{left}}}  ^{0} 
+ h ^0 _{I _{\textup{right}}} ) \,. 
\end{equation*}
\end{proposition}

Each of the five terms on the right are $ L ^{p}$-bounded operators on $ f$, provided $ b\in \textup{BMO}$, so 
that the upper bound on the commutator norm in Theorem~\ref{t.comm} follows as an easy corollary. 
The paraproduct in \eqref{e.oddONe} does not hew 
to our narrow definition of a Paraproduct, but it is degenerate in that it is of signature $ (0,0,0)$, and thus 
even easier to control than the other terms. 

\begin{proof}
Now, $ [ b, \mathfrak H]f= b \mathfrak H f- \mathfrak H (b \cdot f)$.  Apply \eqref{e.product=3para} 
to both of these products.  We see that 
\begin{equation*}
[ b, \mathfrak H]f= 
\sum _{\vec \epsilon = (1,0,0),(0,1,0), (0,0,1)}
\operatorname P ^{\vec \epsilon  } (b , \mathfrak H f) - 
\mathfrak H \operatorname P ^{\vec \epsilon  } (b , f)  \,. 
\end{equation*}
The choices of $ \vec \epsilon = (0,1,0), (0,0,1)$ lead to the first four terms on the right in \eqref{e.px}.  

The terms that require more care are the difference of the two terms in which a $ 1$ falls on a $ b$.  
In fact, we will have 
\begin{equation*}
\operatorname P ^{\vec \epsilon  } (b , \mathfrak H f) - 
\mathfrak H \operatorname P ^{\vec \epsilon  } (b , f) = 
\widetilde {\operatorname  P} ^{0,0,0} (b, f)\,. 
\end{equation*}

To analyze this difference quickly, let us write 
\begin{equation*}
\ip \mathfrak H f, h_I, = \operatorname {sgn} (I) \ip f, h _{\operatorname {Par} (I)}, 
\end{equation*}
where $ \operatorname {Par} (I)$ is the `parent' of $ I$, and $ \operatorname {sgn} (I) = 
1$ if $ I$ is the left-half of $ \operatorname {Par} (I)$, and is otherwise $ -1$. 
This definition follows immediately from the definition of $ g_I$ in \eqref{e.gIdef}.
Now observe that 
\begin{align*}
\ip  {\operatorname P ^{\vec \epsilon  } (b , \mathfrak H f) }, h_I ^{0} , 
&=  
 \ip \mathfrak Hf, {\operatorname P ^{\vec \epsilon  } (b,h_I ^0)} , 
\\
&= \frac {\ip b, h_I ^{1}, } {\sqrt {\lvert  I\rvert }} \cdot \langle \mathfrak H f,  h _{I}^0 \rangle 
\\ 
&= \ip f, h _{\operatorname {Par} (I),  } ^0,  \operatorname {sgn} (I) 
\frac {\ip b, h_I ^{1}, } {\sqrt {\lvert  I\rvert }}  
\end{align*}
And on the other hand, we have 
 \begin{align*}
 \langle \mathfrak H P ^ {1,0,0} (b,f), h_I \rangle &=  
 \frac {\ip b, h _{\operatorname {Par} (I) } ^ {1} , } 
 {\sqrt { \lvert  \operatorname {Par} (I)\rvert } } 
 \operatorname {sgn} (I) \ip f, h ^ {0} _{\operatorname {Par} (I)},    
\end{align*} 
Comparing these two terms, we see that we should examine the term that falls on $ b$.  But a calculation shows that 
\begin{equation*}
\sqrt 2 h _I ^{1} - h _{\operatorname {Par} (I)} ^{1} = - \operatorname {sgn} (I) h _{\operatorname {Par} (I)} ^{0}. 
\end{equation*}
Thus, we see that this difference is of the claimed form.

\end{proof}

\section{The Nehari Theorem} 

We define Hankel operators on the real line. 
On $ L^2(\mathbb R )$, we have the Fourier transform 
\begin{equation*}
\widehat f(\xi )=\int f(x) \operatorname e ^{-i \xi x} \; dx  \,.
\end{equation*}
Define the orthogonal projections onto positive and negative frequencies 
\begin{equation*}
\operatorname P _{\pm} f(x) \eqdef \int _{\mathbb R _{\pm}} \widehat f(\xi ) 
\operatorname e ^{i \xi x} \; dx  \,.
\end{equation*}
Define Hardy spaces $ H ^2   (\mathbb R ) \eqdef \operatorname P _{+} L^2(\mathbb R ) $.
Functions $ f\in H^2(\mathbb R )$ admit an analytic extension to the upper half plane 
$ \mathbb C_+$. As in the case of the disk, it is convenient to refer to functions in 
$ H^2(\mathbb R )$ as \emph{analytic}.

A \emph{Hankel operator with symbol $ b$} is then a linear operator from $ H^2(\mathbb C_+ )$ 
to $ H^2 _{+}(\mathbb C_+ )$ given by $ \operatorname H _{b} \varphi  \eqdef 
\operatorname P_+
\operatorname M_b \overline{\varphi }$. 
This only depends on the analytic part of $ b$.   It is typical to include the notation 
$ \mathbb C _+$ to emphasize the connection with analytic function theory, and the relevant domain 
upon which one is working.  Below, we will suppress this notation.

The  result that we are interested in is: 

\begin{nehari}[\cite {nehari}] \label{t.nehari}
The Hankel operator $ \operatorname H_b$ is bounded from  $ H^2 $ to $ H^2 $ 
iff there is a bounded function $ \beta $ with $ P_+b=P_+\beta $.  Moreover, 
\begin{equation} \label{e.nehari}
\norm \operatorname H_b..=\inf _{\beta \mid \operatorname P_+\beta =\operatorname P_+ b} 
\norm \beta .\infty .
\end{equation}
Less exactly, we have $ \norm \operatorname H_b. . \simeq \norm \operatorname P_+ b . \textup{BMO}.$, 
where we can take the last norm to be non-dyadic $ \textup{BMO}$. 
\end{nehari}

This Theorem was proved in 1954, appealing to the following classical fact. 

\begin{proposition}\label{p.factor}
Each function $ f\in H^1  $ is a product of functions 
$ f_1,f_2\in H^2 $. 
In particular, $ f_1$ and $ f_2$ can be chosen so that 
\begin{equation*}
\norm f.H^1.=\norm f_1.H^2.\norm f_2.H^2. 
\end{equation*}
\end{proposition}

Given a bounded Hankel operator $ \operatorname H _{b}$, we want to show that we 
can construct a bounded function $ \beta $ so that the analytic part of $ b$ and 
$ \beta $ agree. 

This proof is the one found by Nehari \cite {nehari}.  We begin with a basic computation of the 
norm of the Hankel operator $ \operatorname H _{b}$: 
\begin{equation}\label{e.dualNorm}
\begin{split}
\norm \operatorname H _{b}..&=\sup _{\norm \varphi  . H^2 .=1}
\sup _{\norm \psi . H^2 .=1}
\int  \operatorname H_b \psi \cdot \overline{ \varphi} \; dx  
\\
&=\sup _{\norm \varphi  . H^2 .=1}
\sup _{\norm \psi . H^2 .=1} 
\int \operatorname P_+ \operatorname M_b \overline{\psi}  \cdot \overline \varphi \; dx 
\\
&=\sup _{\norm \varphi  . H^2 .=1}
\sup _{\norm \psi . H^2 .=1}\int 
(\operatorname P_+ b)\overline{\psi \cdot \varphi}  \; dx 
\\
&= \sup _{\norm \varphi  . H^2 .=1}
\sup _{\norm \psi . H^2 .=1} \ip (\operatorname P_+ b), \psi \cdot \varphi ,
\end{split}
\end{equation}
But, the $ H^1 =H^2 \cdot H^2 $, 
as we recalled in Proposition~\ref{p.factor}. 
We read from the equality above that 
the analytic part of $ b$ defines a bounded linear functional 
on $ H^1 $ a subspace of $ L^1 $.  

The Hahn Banach Theorem applies, giving us an extension of this linear functional 
to all of $ L^1$, with  the same norm.  But a linear function on $ L^1$ is a bounded 
function, hence we have constructed a bounded function $ \beta $ with the same analytic 
part as $ b$.

The calculation \eqref{e.dualNorm} is more general than what we have indicated here, a point that 
we return to below. 

Let us remark that the $ H ^{p}  $ variant of Nehari's Theorem holds.  On the one hand, 
one has $ H ^{p} \cdot H ^{p'}\subset H ^{1}$, so that the upper bound on the norm $ \norm \operatorname H_b. 
H ^{p} \to H ^{p}.$ follows.  On the other, Proposition~\ref{p.factor} extends to the $ H ^{p}$-$ H ^{p'}$ 
factorization, whence the same argument for the lower bound can be used.

\bigskip

There is a close connection between commutators $ [b, \operatorname H]$ and Hankel operators.  Indeed, 
we have 
\begin{equation}\label{e.[]=Hank}
[b, \operatorname H] = [b, \operatorname H] = 2 \operatorname P_- b \operatorname P_+ - 2\operatorname P_+ b \operatorname
P_- \,. 
\end{equation}
The two terms on the right can be recognized as two  Hankel operators with orthogonal domains and ranges.  
Indeed, keep in mind the elementary identities 
 $ \operatorname P _{+} ^2 = \operatorname P _{+}   $, $ \operatorname P_+ \operatorname P_-=0$, 
 $ \operatorname H=\operatorname I-2 \operatorname P_-$,  and $ [b ,\operatorname I]=0$.   Then, observe 
\begin{align*} 
\operatorname P_+[b, \operatorname H] \operatorname P_- 
&= - 2\operatorname P_+[b, \operatorname P_-] \operatorname P_- 
\\ 
&= -\operatorname P_+ b \operatorname P_- ^2 +  \operatorname P_+\operatorname P_- b \operatorname P_-  
=  -\operatorname P_+ b \operatorname P_-  
\\ 
\operatorname P_- [ b, \operatorname H] \operatorname P_- 
&= \operatorname P_- [ b, \operatorname P_+] \operatorname P_- =0
\end{align*}
There are two additional calculations, which are dual to these and we omit them.

\section{Further Applications} 

The author came to the Haar shift approach  to the commutator
from studies of Multi-Parameter Nehari Theorem \cites{sarahlacey,math.CA/0310348}. 
The paper \cite{math.CA/0601272} surveys 
these two papers. This subject requires an understanding of the structure 
of product $ \textup{BMO}$ that goes beyond the foundational papers of S.-Y.~Chang and R.~Fefferman 
\cites{cf1
,cf2} 
on the subject.   

In particular, as in Nehari's Theorem, the upper bound on the Hankel operator is trivial, as one 
direction of the factorization result is trivial:  $ H ^{2} \cdot H ^{2}\subset H ^{1}$. The lower bound 
is however very far from trivial, as factorization is known to fail in product Hardy spaces.  
Indeed, Nehari's theorem is equivalent to so-called \emph{weak} factorization, one of the points of interest 
in the Theorem.  See \cites{sarahlacey,math.CA/0310348,math.CA/0601272} 
for a discussion of this important obstruction to the proof, and relevant references.

There are different critical ingredients needed for the proof of the lower bound. One of them  is a very precise 
quantitative understanding of the proof of the upper bound.  It is at this point that the techniques 
indicated in this paper are essential.  The fundamentals of the multi-parameter Paraproduct theory were 
developed by Journ\'e \cite{MR88d:42028,MR949001}.  The subject has been revisited recently to develop novel 
Leibnitz rules by Muscalu, Pipher, Tao and Thiele \cites{camil1,camil2}.  Also see \cite{math.CA/0502334}. 

An influential extension of the classical Nehari Theorem to a real-variable setting was found by Coifman, Rochberg 
and Weiss \cite{MR54:843}:  Real-valued $ \textup{BMO} $ on $ \mathbb R ^{n}$ can be characterized in 
terms of commutators with  Riesz Transforms. The real-variable setting implies a complete loss of analyticity, making 
neither bound easy. 
Recently, the author, with Pipher, Petermichl and Wick, have proved the multi-parameter extension of the 
this result \cite{onBeyondCRW}.  This paper includes in it a quantification of the Proposition~\ref{p.[]=} 
to the higher dimensional setting, for (smooth) Calder\'on Zygmund operators $ \operatorname T$: 
$ [b,T]$ is a sum of bounded paraproducts, a crucial Lemma in that paper.  See \cite{onBeyondCRW}*{Proposition 5.11}. 
Such an observation is not new, as it can be found in e.\thinspace g.\thinspace \cite{MR1349230} for instance. 
Still the presentation of Proposition~\ref{p.[]=} in this paper is as simple as any the author is 
aware of in the literature.

%
%
%
%
%

\end{document}